\documentclass[twoside,11pt]{article}

\usepackage{jmlr2e}

\newtheorem{Def}{Definition}[section]

\renewcommand{\c}{\mathbf{c}}

\newcommand{\beqn}{\begin{eqnarray}}
\newcommand{\beqnn}{\begin{eqnarray*}}
\newcommand{\eeqn}{\end{eqnarray}}
\newcommand{\eeqnn}{\end{eqnarray*}}
\newcommand{\GG}{\mathcal{G}}
\def \FF {\mathcal{F}}

\def \FF {\mathcal{F}}
\def \no {\Arrowvert}

\def \R {\mathbb{R}}

\def \R {\mathbb{R}}
\def \E {\mathbb{E}}
\def \N {\mathbb{N}}

\ShortHeadings{Fast rates for Noisy Clustering}{S\'ebastien Loustau}
\firstpageno{1}

\begin{document}

\title{Fast rates for Noisy Clustering}

\author{\name S\'ebastien Loustau \email loustau@math.univ-angers.fr \\
       \addr LAREMA\\
       Universit\'e d'Angers\\
       2 Boulevard Lavoisier,\\
       49045 Angers Cedex, France
       }

\editor{Leslie Pack Kaelbling}

\maketitle

\begin{abstract}
 The effect of errors in variables in empirical minimization is investigated. Given a loss $l$ and a set of decision rules $\mathcal{G}$, we prove a general upper bound for an empirical minimization based on a deconvolution kernel and a noisy sample $Z_i=X_i+\epsilon_i,i=1,\ldots,n$. 
       
      We apply this general upper bound to give the rate of convergence for the expected excess risk in noisy clustering. A recent bound from \citet{levrard} proves that this rate is $\mathcal{O}(1/n)$ in the direct case, under Pollard's regularity assumptions. Here the effect of noisy measurements gives a rate of the form $\mathcal{O}(1/n^{\frac{\gamma}{\gamma+2\beta}})$, where $\gamma$ is the H\"older regularity of the density of $X$ whereas $\beta$ is the degree of illposedness.
\end{abstract}


\section{Introduction}
         Isolate meaningfull groups from the data is an interesting topic in data analysis with applications in many fields, such as biology or social sciences. This unsupervised learning task is known as clustering (see the early work of \citet{hartigan75}). Let $X_1,\ldots ,X_n$ denote i.i.d. random variables with unknown law $P$ on $\R^d$, with density $f$ with respect to a $\sigma$-finite measure $\nu$. The problem of clustering is to assign to each observation a cluster over a finite number of $k$ possible items. From statistical viewpoint, this problem can be endowed into the general and extensively studied problem of empirical minimization (see \citet{vapnik2000}, \citet{kolt}) as follows. Let us consider a class of decision rules $\mathcal{G}$ and a loss function $l:\mathcal{G}\times \R^d$ where $l(g,x)$ measures the loss of $g$ at point $x$. We aim at choosing from the data $X_1,\ldots ,X_n$ a candidate $g\in\GG$ that minimizes the risk functionnal:
\beqn
\label{risk}
R_l(g)=\E l(g,X),
\eeqn
where the expectation is taken over the unknown distribution $P$. For instance the $k$-means algorithm proposes as criterion for partitioning the data the within cluster sum of squares $\mathbf{c}\mapsto \min_{\c}\no x-c_j\no^2$, where in the sequel $\no\cdot\no$ denotes the euclidian norm and $\mathbf{c}=(c_1,\ldots,c_k)$ is the set of possible clusters, with corresponding decision rule $g_\mathbf{c}(x)= \arg\min_{j}\| x-c_j\|$. The performances of a given $g\in\mathcal{G}$ is measured through its non-negative excess risk, given by:
\beqn
\label{excess}
R_l(g)-R_l(g^*),
\eeqn
where $g^*$ is a minimizer over $\GG$ of the risk (\ref{risk}).

A classical way to tackle this issue in the direct case is to consider, if there exists, the Empirical Risk Minimizer (ERM) estimator defined as:
\beqn
\label{erm}
\hat{g}_n=\arg\min_{g\in\GG}R_n(g),
\eeqn
where $R_n(g)$ denotes the empirical risk defined as:
\beqnn
R_n(g)=\frac{1}{n}\sum_{i=1}^nl(g,X_i):=P_nl(g).
\eeqnn
In the sequel the empirical measure of the direct sample $X_1,\ldots , X_n$ will be denoted as $P_n$.  A large litterature (see \citet{vapnik2000} for such a generality) deals with the statistical performances of \ref{erm}) in terms of the excess risk (\ref{excess}). The central point of these papers is to control the complexity of the set $\mathcal{G}$ thanks to VC dimension (\citet{vapnik82}), entropy conditions (\citet{vdg}), or Rademacher complexity assumptions in \citet{localrademacher,kolt} ( see also \citet{nedelec,svm} in supervised classification). The main probabilistic tool for this problem is the statement of uniform concentration of the empirical measure to the true measure. This can be easily seen using the so-called Vapnik's bound:
\beqn
\label{vb}
R_l(\hat{g}_n)-R_l(g^*)&\leq& R_l(\hat{g}_n)-R_n(\hat{g}_n)+R_n(g^*)-R_l(g^*)\noindent\\
&\leq & 2\sup_{g\in\GG}|(R_n-R_l)(g)|=2\sup_{g\in\GG}|(P_n-P)l(g)|.
\eeqn
It is important to note that (\ref{vb}) can be improved using a local approach (see \citet{toulouse}) which consists in reducing the supremum to a neighborhood of $g^*$. We do not develop this important refinement in this introduction for the sake of concision whereas it is the main ingredient of the literature cited above. It allows to get fast rates of convergence.\\
 
In this paper the framework is essentially different since we observe a corrupted sample $Z_1,\ldots,Z_n$ such that:
\beqn
\label{noisy}
Z_i=X_i+\epsilon_i,\,i=1,\ldots n,
\eeqn 
where the $\epsilon_i$'s are i.i.d. $\R^d$-random variables with density $\eta$ with respect to the Lebesgue measure. As a result, from (\ref{noisy}, the empirical measure $P_n=\frac{1}{n}\sum_{i=1}^n\delta_{X_i}$ is unobservable and standard ERM (\ref{erm}) is not available. Unfortunately, using the corrupted sample $Z_1,\ldots,Z_n$ in standard ERM (\ref{erm}) seems problematic since:
\beqnn
\tilde{P}_nl(g):=\frac{1}{n}\sum_{i=1}^nl(g,Z_i)\longrightarrow \E l(g,Z)\not= R_l(g).
\eeqnn
Due to the action of the convolution operator, the empirical measure of the indirect sample, defined as $\tilde{P}_n:=\frac{1}{n}\sum_{i=1}^n\delta_{Z_i}$ differs from $P_n$ and we are faced to an ill-posed inverse problem. Note that this problem has been recently considered in \citet{pinkfloyds} in discriminant analysis and in a more general supervised statistical learning context in \citet{noisykolt}. The main idea to get optimal upper bounds is to consider an empirical risk based on kernel deconvolution estimators.

In this paper, we propose to adopt a comparable strategy in unsupervised statistical learning. To this end, we propose to construct a kernel deconvolution estimator of the density $f$ of the form:
\beqn
\label{dke}
\hat{f}_n(x)=\frac{1}{n}\sum_{i=1}^n\frac{1}{\lambda}\mathcal{K}_\eta\left(\frac{Z_i-x}{\lambda}\right),
\eeqn
where $\mathcal{K}_\eta$ is a deconvolution kernel and $\lambda$ is a regularization parameter (see Section \ref{sec2} for details). 
Given this estimator, we construct an empirical risk by plugging (\ref{dke}) into the true risk $R_l(g)$ to get
a so-called deconvolution empirical risk minimization given by:
\beqn
\label{lerm}
\arg\min_{g\in\mathcal{G}}R_n^\lambda(g)\mbox{ where }R_n^\lambda(g):=\frac{1}{n}\sum_{i=1}^nl_\lambda(g,Z_i),
\eeqn
whereas $l_\lambda(g,z)$ is a convolution of the loss function $l(g,\cdot)$ given by: 
$$
l_\lambda(g,z)=\int\frac{1}{\lambda}\mathcal{K}_\eta\left(\frac{z-x}{\lambda}\right)l(g,x)\nu(dx).
$$ 
Note that in case no such minimum exists, we can consider $\delta$-approximate minimizers as in \citet{empimini}. 

In order to study the performances of a solution of (\ref{lerm}, it is possible to use the empirical process machinery in the spirit of \citet{kolt,empimini,svm}. In the presence of indirect observations, for $\hat{g}_n^\lambda$ a solution of the minimization of (\ref{lerm}, we have:
\beqn
\label{nvb}
R_l(\hat{g}_n^\lambda)-R_l(g^*)&\leq& R_l(\hat{g}_n^\lambda)-R_n^\lambda(\hat{g}_n^\lambda)+R_n^\lambda(g^*)-R_l(g^*)\noindent\\
&\leq& R_l^\lambda(\hat{g}_n^\lambda)-R_n^\lambda(\hat{g}_n^\lambda)+ R_n^\lambda(g^*)-R_l^\lambda(g^*)+ (R_l-R_l^\lambda)(\hat{g}_n^\lambda -g^*)\noindent\\
&\leq & \sup_{g\in\GG}|(R_n^\lambda-R_l^\lambda)(g^*-g)|+\sup_{g\in\GG}|(R_l^\lambda-R_l(g-g^*)|,
\eeqn
where in the sequel, under integrability conditions and using Fubini: 
\beqn
\label{rl}
R_l^\lambda(g)=\E R_n^\lambda(g)= \int l(g,x)\E\frac{1}{\lambda} \mathcal{K}_\eta\left(\frac{Z-x}{\lambda}\right)\nu(dx).
\eeqn
Bounds (\ref{nvb}) are comparable to (\ref{vb}) for the direct case. There consist in two terms: 
\begin{itemize}
\item A variance term $\sup_{g\in\GG}|(R_n^\lambda-R_l^\lambda)(g^*-g)|$ related to the estimation of $R_l(g)$ using an empirical couterpart. This term will be controled using standard tools from empirical process theory, namely a local approach in the spirit of \citet{kolt}. Here the empirical process is indexed by a class of functions depending on a smoothing parameter.
\item A bias term $\sup_{g\in\GG}|(R_l^\lambda-R_l)(g-g^*)|$ due to the estimation procedure using kernel deconvolution estimator. It seems to be related to the usual bias term in nonparametric density deconvolution since we can see coarselly that:
\beqnn
R_l^\lambda(g)-R_l(g)=\int l(g,x)\left[\E\frac{1}{\lambda} \mathcal{K}_\eta\left(\frac{Z-x}{\lambda}\right)-f(x)\right]\nu(dx).
\eeqnn
\end{itemize}
The choice of $\lambda$ is crucial in the decomposition (\ref{nvb}. We will show below that the variance term grows when $\lambda$ tends to zero whereas the bias term vanishes. Parameter $\lambda$ has to be chosen as a trade-off between these two terms, and as a consequence will depend on unknown parameters. The problem of adaptation is not adressed in this paper but is an interesting future direction.\\

The paper is organized as follows. In Section 2, we propose to give a general upper bound for (\ref{lerm}, generalizing the results of \citet{kolt} to indirect observations. Note that all the material of Section \ref{sec2} is largely inspired from \citet{noisykolt} and gives an unsupervised counterpart of the previous results. Section 3 gives a direct application of the result of Section 2 in clustering by giving rates of convergence for a new deconvolution $k$-means algorithm. Fast rates of convergence are proposed which generalize the recent fast rates proposed in \citet{levrard} in the direct case. 

\section{General Upper bound}
\label{sec2}
In this section we propose an upper bound for the expected excess risk of the estimator:
\beqn
\label{decerm}
\hat{g}_n^\lambda:=\arg\min \frac{1}{n}\sum_{i=1}^nl_\lambda(g,Z_i),
\eeqn
where $l_\lambda(g,z)$ is construct as follows.\\
Let us introduce $\mathcal{K}=\prod_{i=1}^d \mathcal{K}_j:\R^d \to \R$ a $d$-dimensional function defined as the product of $d$ unidimensional function $\mathcal{K}_j$. Then if we denote by $\lambda=(\lambda_1,\dots,\lambda_d)$ a set of (positive) bandwidths and by $\FF[\cdot]$ the Fourier transform, we define $\mathcal{K}_\eta$ as:
\begin{eqnarray}
\mathcal{K}_{\eta} & : & \R^d \to \R \nonumber \\
& & t \mapsto \mathcal{K}_\eta(t) = \FF^{-1}\left[ \frac{\FF[\mathcal{K}](\cdot)}{\FF[\eta](\cdot/\lambda)}\right](t).
\label{dk}
\end{eqnarray}
Moreover in the sequel we restrict the study to a compact set $K\subset\R^d$ and define $l_\lambda(g,z)$ as
$$
l_\lambda(g,z)=\int_K\frac{1}{\lambda}\mathcal{K}_\eta\left(\frac{z-x}{\lambda}\right)l(g,x)\nu(dx),
$$ 
where we write with a slight abuse of notation $\frac{1}{\lambda}\mathcal{K}_\eta\left(\frac{z-x}{\lambda}\right)$ for $\frac{1}{\Pi_{i=1}^d\lambda_i}\mathcal{K}_\eta\left(\frac{z_1-x_1}{\lambda_1},\ldots,\frac{z_d-x_d}{\lambda_d}\right)$. The restriction to a compact $K$ allows to control the variance in decomposition (\ref{nvb}) thanks to Lemma \ref{lip} (we refer the reader to \citet{noisykolt} for a discussion).

Finally in the sequel for the sake of simplicity we restrict ourselves to moderately ill-posed inverse problem and introduce the following assumption:
 \\

\noindent
\textbf{Noise Assumption}: There exist $(\beta_1,\dots,\beta_d)'\in \R_+^d$ such that for all $i\in \lbrace 1,\dots, d \rbrace$,
$$ \left| \mathcal{F}[\eta_i](t) \right| \sim |t|^{-\beta_i},  \mathrm{as} \ t\to +\infty.$$
Moreover, we assume that $\mathcal{F}[\eta_i](t) \not = 0$ for all $t\in \R$ and $i\in \lbrace 1,\dots, d \rbrace$.
\\

\hspace{-0.8cm} Assumption \textbf{(NA)} deals with the asymptotic behavior of the characteristic function of the noise distribution. These kind of restrictions are standard in deconvolution problems (see \citet{Fan,meister,butucea}). Note that straightforward modifications allow to consider severely ill-posed inverse problems, where the asymptotic behavior of the characteristic function of $\epsilon$ decreases exponentially to zero.\\

Under \textbf{(NA)}, the goal is to control the two terms of (\ref{nvb}, namely the bias term and the variance term. The variance term is reduced to the study of the increments of the empirical process:
\beqnn
\nu_n^\lambda(g)=\frac{1}{\sqrt{n}}\sum_{i=1}^n l_\lambda(g,Z_i)-\E l_\lambda(g,Z).
\eeqnn 
It will be controled thanks to a version of Talagrand's inequality due to \citet{bousquet}. However it is important to note that here the empirical process is indexed by the class of functions $\{z\mapsto l_\lambda(g,z),g\in\GG\}$, which depends on a regularization parameter $\lambda\in\R^d_+$. This parameter will be calibrated as a function of $n$ so  Talagrand's type inequality has to be used in a careful way. For this purpose, we need the following lemma.
\begin{lemma}
\label{lip}
Suppose $f\geq c_0>0$ on $K$. Then if \textbf{(NA)} holds, and $\mathcal{K}$ has compactly supported Fourier transform, we have:
\begin{description}
\item[(i)] $l(g)\mapsto l_\lambda(g)$ is Lipschitz with respect to $\lambda$:
\beqnn
\exists C_1>0:\forall g,g'\in\GG,\,\Arrowvert l_\lambda(g)-l_\lambda(g')\Arrowvert_{L_2(\tilde{P})}\leq C_1 \Pi_{i=1}^d\lambda_i^{-\beta_i}\Arrowvert l(g)-l(g')\Arrowvert_{L_2(P)}.
\eeqnn
\item[(ii)] $\{l_\lambda(g),g\in\GG\}$ is uniformly bounded:
$$
\exists C_2>0:\sup_{g\in\GG}\| l_\lambda(g,\cdot)\|_\infty\leq C_2 \Pi_{i=1}^d\lambda_i^{-\beta_i-1/2}.
$$
\end{description}
\end{lemma}
The proof of this result is presented in \citet{noisykolt} in a slightly different framework. Note that the assumption on $f$ to be strictly positive on $K$ appears for some technical reasons in the proof and could be avoided in some cases (see the discussion in \citet{noisykolt}).\\The Lipschitz property \textbf{(i)} is a key ingredient to control the complexity of the class of functions $\{l_\lambda(g),g\in\GG\}$ thanks to standard complexity arguments applied to the loss class $\{l(g),g\in\GG\}$. Finally \textbf{(ii)} is necessary to apply Talagrand's type inequality to the empirical process $g\mapsto\nu_n^\lambda(g)$ above. 

To control the excess risk of the procedure, we also need to control the bias term defined in (\ref{nvb}) thanks to Lemma \ref{bias}) below.
\begin{lemma}
\label{bias}
Suppose $f\in\Sigma(\gamma,L)$ the H\"older class of $\lfloor\gamma\rfloor$-fold continuously differentiable functions on $R^d$ satisfying the H\"older condition. Let $\mathcal{K}$ a kernel of order $\lfloor\gamma\rfloor$ with respect to $\nu$. Then if $l(g,\cdot)\in L_1(\nu,\R^d)$, we have:
\beqnn
\forall g,g'\in\GG,\,\left|(R_l-R^\lambda_l)(g-g')\right|\leq C\sum_{i=1}^d\lambda_i^\gamma.
\eeqnn
\end{lemma}
The proof is presented in \citet{noisykolt} and is omitted. Finally to state fast rates, we also require an additional assumption over the distribution $P$. 
\begin{Def}
\label{KKb}
We say that $\mathcal{F}$ is a Bernstein class with respect to $P$ if there exists $\kappa_0\geq 0$ such that for every $f\in\FF$:
\beqnn
\no f\no_{L_2(P)}^2\leq \kappa_0[\E_P f].
\eeqnn
\end{Def}
This notion of Bernstein class first appears in \citet{empimini} in a more general form. Definition \ref{KKb} corresponds to the ideal case where $\kappa=1$. This assumption can be related to the well-known margin assumption in supervised classification, introduced in \citet{mammen}. Section \ref{application} proposes an unsupervised version of this hypothesis.\\
From technical viewpoint, this requirement arises naturally in the proof when we want to apply functional Bernstein's inequality such as Talagrand's inequality. If we consider the loss class $\FF=\{l(g)-l(g^*),g\in\mathcal{G}\}$, Definition \ref{KKb} gives a perfect variance-risk correspondance.

We are now on time to state the main result of this section.
\begin{theorem}
\label{mainresult}
Suppose \textbf{(NA)} holds and assumptions of Lemma 1-2 hold. Suppose $\{l(g)-l(g^*),g\in\mathcal{G}\}$ is Bernstein w.r.t. $P$ where $g^*\in\arg\min_\mathcal{G}R_l(g)$ is unique and there exists $0<\rho<1$ such that for every $\delta>0$:
\beqn
\label{modulus}
\E\sup_{g,g'\in\mathcal{G}(\delta)}\left|(\tilde{P}-\tilde{P}_n)(l_\lambda(g)-l_\lambda(g'))\right|\lesssim \frac{\Pi_{i=1}^d\lambda_i^{-\beta_i}}{\sqrt{n}}\delta^{\frac{1-\rho}{2}},
\eeqn
where $\mathcal{G}(\delta)=\{g\in\mathcal{G}:R_l(g)-R_l(g^*)\leq \delta\}$.\\
Then estimator $\hat{g}=\hat{g}_n^\lambda$ defined in (\ref{decerm}) is such that:
\beqnn 
\E R_l(\hat{g})-R_l(g^*)\leq Cn^{-\frac{\gamma}{\gamma(1+\rho)+2\bar{\beta}}},
\eeqnn
where $\bar{\beta}=\sum_{i=1}^d\beta_i$ and $\lambda=(\lambda_1,\ldots,\lambda_d)$ is chosen as:
$$
\lambda_i\sim n^{-\frac{1}{\gamma(1+\rho)+2\bar{\beta}}},\forall i=1,\ldots d.
$$.
\end{theorem}
The proof of this result iterates a version of Talagrand's inequality due to \citet{bousquet}. It is presented in Section 5. coarsely speaking, Lemma \ref{lip}, gathering with the complexity assumption (\ref{modulus}, leads to a control of the variance term in decomposition (\ref{nvb}. Then Lemma \ref{bias} gives the order of the bias term. The choice of $\lambda$ explicited in Theorem \ref{mainresult} trades off these two terms, and gives the excess risk bound.

Note that the rates of convergence in Theorem \ref{mainresult} generalize previous results. When $\epsilon=0$ Theorem \ref{mainresult} gives fast rates of convergence between $\mathcal{O}(1/\sqrt{n})$ to $\mathcal{O}(1/n)$ depending on the complexity parameter $\rho>0$ in (\ref{modulus}, which can be related with entropy or Rademacher complexities of the hypothesis set $\mathcal{G}$. The effect of the inverse problem depends on the asymptotic behavior of the characteristic function of the noise distribution $\epsilon$. This is rather standard in the statistical inverse problem literature (\citet{Fan} or \citet{meister}). 

Moreover the control of the modulus of continuity in (\ref{modulus}) is specific to the indirect framework and depends on the smoothing parameter $\lambda$. A comparable hypothesis arises in the direct case in \citet{kolt}, up to the constant depending on $\lambda$. It appears that it will be satisfied in our application using standard statistical learning argues, such as maximal inequalities and chaining.

Finally note that the complexity parameter involved in assumption (\ref{modulus}) is smaller than the complexity proposed in \citet{kolt} or \citet{noisykolt}. Here the supremum is taken over the set $\{g,g'\in\mathcal{G}(\delta)\}\subset \{g,g'\in\mathcal{G}:P(l(g)-l(g'))^2\leq c\delta\}$ provided that $\{l(g)-l(g^*),g\in\mathcal{G}\}$
 is Bernstein with respect to $P$ according to Definition \ref{KKb}. This indexing set is related to the localization's technique used in Theorem \ref{mainresult}, namely a localization based on the excess risk instead of the $L_2(P)$-norm. This refinement is necessary to derive fast rates in Section \ref{application} (see \citet{empimini} for a related discussion).   
\section{Application to noisy clustering}
\label{application}
Clustering is a basic problem in statistical learning where independent random variables $X_1,\dots,X_n$ are observed, with common source distribution $P$. The aim is to construct clusters to classify these data. However in many real-life situations, direct data $X_1,\ldots ,X_n$ are not available and measurement errors occur. Then we observe a corrupted sample $Z_i=X_i+\epsilon_i,i=1,\ldots n$ with unknown noisy distribution $\tilde{P}$. The problem of noisy clustering is to learn clusters for the direct dataset $X_1,\ldots, X_n$ when only a contaminated version $Z_1,\ldots ,Z_n$ is observed.

To frame the noisy clustering problem as a statistical learning one, we first introduce the following notation.
         Let $\mathbf{c}= (c_1, \ldots, c_k)\in\mathcal{C}$ the set of possible clusters, where $\mathcal{C}\subseteq\R^{dk}$ and $X\in\R^d$. The loss function $\gamma:\R^{dk}\times \R^d$ is defined as:
$$
\gamma(\mathbf{c},x)=\min_{j=1,\ldots k}\no x-c_j\no^2,
$$
and the corresponding true risk or clustering risk is $R(\mathbf{c})=\E\gamma(\mathbf{c},X)$. The performances of the empirical minimizer $\hat{\mathbf{c}}_n=\arg\min_\mathcal{C}P_n\gamma(\mathbf{c})$ (also called $k$-means clustering algorithm) have been widely studied in the literature. Consistency was shown by \citet{pollard81} when $\E\no X\no^2< \infty$ whereas \citet{llz} or \citet{biau} gives rates of convergence of the form $\mathcal{O}(1/\sqrt{n})$ for the excess clustering risk defined as $R(\hat{\mathbf{c}}_n)-R(c^*)$, where $c^*\in\mathcal{M}$ the set of all possible optimal clusters. More recently, \citet{levrard} proposes fast rates of the form $\mathcal{O}(1/n)$ under Pollard's regularity assumptions. It improves a previous result of \citet{gg}. The main ingredient of the proof is a localization argument in the spirit of \citet{svm}. 

In this section, we study the problem of clustering where we have at our disposal a corrupted sample $Z_i=X_i+\epsilon_i$, $i=1,\ldots ,n$ where the $\epsilon_i$'s are i.i.d. with density $\eta$ satisfying \textbf{(NA)}. For this purpose, we introduce the following deconvolution empirical minimization:
\beqn
\label{noisykmeans}
\arg\min_{\mathbf{c}\in\mathcal{C}}\frac{1}{n}\sum_{i=1}^n\gamma_\lambda(\mathbf{c},Z_i),
\eeqn
where $\gamma_\lambda(\mathbf{c},z)$ is a deconvolution cluster sum of squares defined as:
\beqnn
\gamma_\lambda(\mathbf{c},z)=\int_{K}\frac{1}{\lambda}\mathcal{K}_\eta\left(\frac{z-x}{\lambda}\right)\min_{j=1,\ldots k}\|x-c_j\|^2dx,
\eeqnn
for $\mathcal{K}_\eta$ the deconvolution kernel of Section \ref{sec2} and $\lambda=(\lambda_1,\ldots,\lambda_d)\in\R^d_+$ a set of positive bandwidths chosen later on. Note that here the existence of a minimizer in (\ref{noisykmeans}) could be managed as in \citet{existence} for the direct case. We investigate the generalization ability of the solution of (\ref{noisykmeans}) in the context of Pollard's regularity assumptions, thanks to the noisy empirical minimization results of Section \ref{sec2}. To this end, we will use the following assumptions on the source distribution $P$.
                   \\

\noindent
\textbf{Assumption 1 (Boundedness assumption)}: The distribution $P$ is such that:
$$P({\mathcal{B}}(0,M)) = 1,$$
where  $\mathcal{B}(0,M)$ denote the closed ball of radius $M$, with $M \geq 0$.
\\

\hspace{-0.8cm} Note that \textbf{(A1)} imposes a boundedness condition on the random variable $X$. We will also need the following regularity requirement, first introduced in \citet{Pollard82}.
                    \\

\noindent
\textbf{Assumption 2 (Pollard's regularity condition)}: The distribution $P$ satisfies the following two conditions:
                  
                  \begin{enumerate}
                  \item $P$ has a continuous density $f$ with respect to Lebesgue measure on $\mathbb{R}^d$,
                  \item The Hessian matrix of $ \mathbf{c} \longmapsto P \gamma (\mathbf{c},.) $ is positive definite for all optimal vector of clusters $\mathbf{c}^*$.
                  \end{enumerate}
                  It is easy to see that using the compactness of $\mathcal{B}(0,M)$, \textbf{(A1)-(A2)} ensures that there exists only a finite number of optimal clusters $\c^*\in\mathcal{M}$. This number is denoted as $|\mathcal{M}|$ in the rest of the paper.\\
Moreover, Pollard's conditions can be interpreted as follows. Denote $\partial V_i$ the boundary of the Voronoi cell $V_i$ associated with $c_i$, for $i=1, \ldots, k$. Then \citet{levrard} has shown that a sufficient condition to have \textbf{(A2)} is to control the sup-norm of $f$ on the union of all possible $|\mathcal{M}|$ boundaries $\partial V^{*,m}=\cup_{i=1}^k\partial V^{*,m}_i$, associated with $c^*_m\in\mathcal{M}$ as follows:
$$
\|f_{|\cup_{m=1}^\mathcal{M} \partial V^{*,m}}\|_\infty\leq c(d)M^{d+1}\inf_{m=1,\ldots,|\mathcal{M}|,i=1,\ldots k}P(V_i^{*,m}),
$$
where $c(d)$ is a constant depending on the dimension $d$. As a result, \textbf{(A2)} is guaranteed when the source distribution $P$ is well concentrated around its optimal clusters, which is related to well-separated classes. From this point of view, Pollard's regularity conditions can be related to the margin assumption in binary classification (see \citet{tsybakov2004}). We have in fact the following lemma due to \citet{gg}.
\begin{lemma}[\citet{gg}]
\label{lemap}
Suppose \textbf{(A1)-(A2)} are satisfied. Then, for any $\mathbf{c}\in\mathcal{B}(0,M)$:
$$
\mathrm{var}\left( \gamma(\mathbf{c},\cdot)-\gamma(\mathbf{c}^*(\c),\cdot)\right)\leq C_1\|\c-\c^*(\c)\|^2\leq C_1C_2\left(R(\c)-R(\c^*(\c))\right),
$$
where $c^*(\c)\in\arg\min_{\c^*}\|\c-\c^*\|$.
\end{lemma}                
Lemma \ref{lemap} is useful to derive fast rates of convergence for two reasons.
 
Firstly, if we compil these two inequalities, we get a control of the variance of the excess loss $\gamma(\mathbf{c})-\gamma(\c^*(\c))$ thanks to the excess clustering risk $R(\c)-R(\c^*(\c))$. Note that if $R(\c)-R(\c^*(\c))\leq 1$, it is clear that the loss class  $\{\gamma(\mathbf{c})-\gamma(\c^*(\c)),c\in\mathcal{C}\}$ is Bernstein according to Definition \ref{KKb} since we have coarsely:
$$
P(\gamma(\mathbf{c},\cdot)-\gamma(\c^*(\c),\cdot))^2\leq \mathrm{var} \left(\gamma(\mathbf{c},\cdot)-\gamma(\c^*(\c),\cdot)\right)+(R(\c)-R(\c^*(\c)))^2\leq (C_1C_2+1)\left(R(\c)-R(\c^*(\c))\right).
$$

Moreover the second inequality of Lemma \ref{lemap} is necessary to control the complexity involved in Section \ref{sec2} thanks to the following lemma:
\begin{lemma}
\label{complexity}
Suppose \textbf{(A1)-(A2)} are satisfied. Suppose $\E\|\epsilon\|^2<\infty$. Then: 
$$
\mathbb{E} \sup_{(\c,\c^*)\in\mathcal{C}\times\mathcal{M}, \|\mathbf{c}-\mathbf{c}^*\|^2 \leq \delta} {|\tilde{P}_n-\tilde{P}|(\gamma_\lambda(\mathbf{c}^*,.) - \gamma_\lambda(\mathbf{c},.))} \leq C \Pi_{i=1}^d\lambda_i^{-\beta_i}\frac{\sqrt{\delta}}{\sqrt{n}},
$$
where $C$ is a positive constant depending on $M,k,d,\tilde{P},K,\eta$.
\end{lemma} 
Note that $\E\|\epsilon\|^2<\infty$ comes from \citet{Pollard82}. Gathering with \textbf{(A1)}, it gives $\E\|Z\|^2<\infty$ and allows to deal with indirect observations. The proof of Lemma \ref{complexity} is presented in Section \ref{proofs}. It is based on \citet{Pollard82} extended to the noisy setting. Under \textbf{(A2)} and provided that $\E\|\epsilon\|^2<\infty$, we use the following approximation of the convolution loss function $\gamma_\lambda(\cdot,x)$ at any point $\mathbf{c}\in\mathcal{C}$:
\beqn
\label{taylor}
\gamma_\lambda(\mathbf{c},z)=\gamma_\lambda(\mathbf{c}^*,z)+\left\langle \mathbf{c}-\mathbf{c}^*,\nabla_{\mathbf{c}}\gamma_\lambda(\mathbf{c}^*,z) \right\rangle + \|\mathbf{c}-\mathbf{c}^*\|R_\lambda(\mathbf{c}^*,\mathbf{c}-\mathbf{c}^*,z),
\eeqn
where $\nabla_{\mathbf{c}}\gamma_\lambda(\mathbf{c}^*,z)$ is the gradient of $\mathbf{c}\mapsto\gamma_\lambda(\mathbf{c},z)$ at point $\c^*$ and $R_\lambda(\mathbf{c}^*,\mathbf{c}-\mathbf{c}^*,z)$ is a residual term (see \citet{Pollard82} for details). With (\ref{taylor}, the complexity term is controled with a maximal inequality due to \citet{saintflour}, gathering with a chaining method.

We are now on time to state the main result of this section.
\begin{theorem}
\label{thm:application}
Assume \textbf{(NA)} holds, $P$ satisfies \textbf{(A1)-(A2)} with density $f\in\Sigma(\gamma,L)$ and $\E\|\epsilon\|^2<\infty$. Then, denoting by $\hat{\mathbf{c}}^\lambda_n$ a solution of (\ref{noisykmeans}, we have, for any $c^*\in\mathcal{M}$:
                  
                  \[
                  \mathbb{E} R(\hat{\mathbf{c}}^\lambda_n)-R(\c^*) \leq Cn^{-\frac{\gamma}{\gamma+2\bar{\beta}}},
                  \]
                  where $\bar{\beta}=\sum_{i=1}^d\beta_i$, $C$ is a positive constant whereas $\lambda=(\lambda_1,\ldots,\lambda_d)$ is chosen as:
                 $$
\lambda_i\sim n^{-\frac{1}{\gamma+2\bar{\beta}}},\forall i=1,\ldots d.
$$
                  \end{theorem}
The proof is a direct application of Section \ref{sec2} when $|\mathcal{M}|=1$ whereas when $|\mathcal{M}|\geq 2$, a more sophisticated geometry has to be considered (see Section \ref{proofs} for details). Some remarks are in order.\\
Rates of convergence of Theorem \ref{thm:application} are fast rates when $2\bar{\beta}<\gamma$. It generalizes the result of \citet{levrard} to the errors-in-variables case since we can see coarsely that rates to the order $\mathcal{O}(1/n)$ are reached when $\epsilon=0$. Here the price to pay for the inverse problem is the quantity $2\sum_{i=1}^d \beta_i$, related to the tail behavior of the characteristic function of the noise distribution $\eta$ in \textbf{(NA)}. This rate corresponds to the ideal case where $\rho=0$ in Section \ref{sec2}, due to the finite-dimensional structure of the set of clusters $\mathcal{C}=\{\c=(c_1,\ldots,c_k),c_i\in\R^d\}$. An interesting extension is to consider richer classes such as kernel classes (see \citet{mendelsonkernel}) and to deal with kernel $k$-means.\\
Lower bounds of the form $\mathcal{O}(1/\sqrt{n})$ have been stated in the direct case by \citet{lbclustering} for general distribution. An open problem is to derive optimality of Theorem \ref{thm:application}, even in the direct case where $\epsilon=0$. For this purpose, we need to construct configurations where both Pollard's regularity assumption and noise assumption \textbf{(NA)} could be used in a careful way. In this direction \citet{pinkfloyds} proposes lower bounds in a supervised framework under both margin assumption and \textbf{(NA)}.

\section{Conclusion}
This paper can be seen as a first attempt into the study of both empirical minimization and clustering with errors-in-variables. Many problems could be considered in future works, from theoretical or practical point of view.

In the problem of empirical minimization with errors-in-variables, we provide the order of the expected excess risk, depending on the complexity of the hypothesis space, the regularity of the direct observations and the degree of ill-posedness. For the sake of concision, Theorem \ref{mainresult} only consider particular Bernstein classes and empirical minimization based on a deconvolution kernel estimator. A higher level of generality can be derived from \citet{noisykolt} but is out of the scope of the present paper. 

The performances of our deconvolution $k$-means algorithm is obtained thanks to a localization principle due to \citet{kolt}, where proofs iterate a Talagrand's inequality due to \citet{bousquet}. With such a study, \citet{kolt} provides the order of the excess risk in the direct case and allows to recover most of the recent results in the statistical learning context. There is nice hope that many statistical learning problem when dealing with indirect observations could be solved with similar argues. 

In the problem of noisy $k$-means clustering, we propose fast rates of convergence to the order of $\mathcal{O}(1/n^{\frac{\gamma}{\gamma+2\bar{\beta}}})$. Theorem \ref{thm:application} is a direct application of the result of Theorem \ref{mainresult} to the problem of clustering with Pollard's regularity assumptions and bounded source. It generalizes a recent result of \citet{levrard} where fast rates are stated for direct observations. 

From practical viewpoint, this work proposes an empirical minimization to deal with the problem of noisy clustering. However the procedure (\ref{noisykmeans}) is not adaptive in many sense. Of course the dependency on the noise distribution $\eta$ can be explored in a future work, and can be associated with the problem of unknown operator in the statistical inverse problem literature (see \citet{marteau1} or \citet{cavalierhengartner}). Moreover the empirical minimization depends on the H\"older regularity of the density of the source distribution through the choice of the bandwidths $\lambda_i,i=1,\ldots d$. However for practical experiments, any data-dependent model selection procedure can be performed, such as cross-validation.            
\section{Proofs}
In all the proofs constant $C>0$ may vary from line to line.
\label{proofs}
\subsection{Proof of Theorem \ref{mainresult}}
The proof uses the following intermediate lemma. 
\begin{lemma}
\label{interm}
Suppose $\{l_\lambda(g),g\in\mathcal{G}\}$ is such that $\sup\no l_\lambda(g)\no_\infty\leq K(\lambda)$. Define:
\beqnn
U_{n}^\lambda(\delta_j,t):=K\left[\phi_n^\lambda(\GG,\delta_j)+\sqrt{\frac{t}{n}}D^\lambda(\GG,\delta_j)+\sqrt{\frac{t}{n}(1+K(\lambda))\phi_n^\lambda(\GG,\delta_j)}+\frac{t}{n}\right],
\eeqnn
\beqnn
\phi_n^\lambda(\GG,\delta_j):=\E \sup_{g,g'\in\GG (\delta_j)}|\tilde{P}_n-\tilde{P}|[l_\lambda(g)-l_\lambda(g')],
\eeqnn
\beqnn
D^\lambda(\GG,\delta_j):=\sup_{g,g'\in\GG (\delta_j)}\sqrt{\tilde{P}(l_\lambda(g)-l_\lambda(g'))^2},
\eeqnn
where $\delta_j=q^{-j},j\in\N^*$, for some $q>0$.\\
Then $\forall \delta\geq \delta_n^\lambda(t)$, we have for $\hat{g}=\hat{g}_n^\lambda$ defined in (\ref{decerm}:
\beqnn
\mathbb{P}(R_l(\hat{g})-R_l(g^*)\geq \delta)\leq c(\delta,q)e^{-t},
\eeqnn
where:
$$
\delta_n^\lambda(t)=\left(\inf\left\{\delta>0:\sup_{\delta_j\geq \delta}\frac{U_n(\delta_j,t)}{\delta_j}\leq\frac{1}{2q}\right\}\right)\vee \left(4q\sup_{g,g'\in\mathcal{G}}(R_l-R_l^\lambda)(g-g')\right).
$$
\end{lemma}
The proof is a straightforward modification of the proof of Lemma 2 in \citet{noisykolt}.
\begin{proof}[of Theorem \ref{mainresult}]
First note that, in dimension $d=1$ for simplicity:
\beqnn
U_n^\lambda(\delta,t)\leq C\left(\phi_n^\lambda(\GG,\delta)+\sqrt{\frac{t}{n}\phi_n^\lambda(\GG,\delta)(1+\lambda^{-\beta-1/2 })}+\sqrt{\frac{t}{n}}D^\lambda(\GG,\delta)+\frac{t}{n}\right).
\eeqnn
Using Definition \ref{KKb}, gathering with the complexity assumption over $\tilde{\omega}_n(\mathcal{G},\delta)$, we have:
\beqnn
\phi_n^\lambda(\GG,\delta)\leq \E \sup_{g,g'\in\GG (\delta)}|\tilde{P}_n-\tilde{P}|[l_\lambda(g)-l_\lambda(g')]\leq C\frac{\lambda^{-\beta}}{\sqrt{n}}\delta^{\frac{1-\rho}{2}}.
\eeqnn
A control of $D^\lambda(\GG,\delta)$ using Lemma \ref{lip}, gathering with Definition \ref{KKb} leads to:
\beqnn
U_n^\lambda(\delta,t)\leq C\left(\frac{\lambda^{-\beta}}{\sqrt{n}}\delta^{\frac{1-\rho}{2}}+\frac{\lambda^{-\beta/2}}{n^{3/4}}\delta^{\frac{1-\rho}{4}}\sqrt{\lambda^{-\beta-1/2}t}+\sqrt{\frac{t}{n}} \lambda^{-\beta}\delta^{\frac{1}{2}}+\frac{t}{n}\right).
\eeqnn
We hence have from an easy calculation:
\beqnn
\delta_n^\lambda(t)\leq C\left(\frac{\lambda^{-\beta}}{\sqrt{n}}\right)^{\frac{2}{1+\rho}}.
\eeqnn
To get the result we apply Lemma \ref{interm} with:
$$
\delta=C(1+t)\left(\frac{\lambda^{-\beta}}{\sqrt{n}}\right)^{\frac{2}{1+\rho}},
$$
noting that the choice of $\lambda$ warrants that:
$$
\lambda^{\gamma}\leq C(1+t)\left(\frac{\lambda^{-\beta}}{\sqrt{n}}\right)^{\frac{2}{1+\rho}}.
$$
Same arguments conclude the proof for $d\geq 2$.
\end{proof}               
\subsection{Proof of Lemma \ref{complexity}}

The proof follows \citet{levrard} applied to the noisy setting. First note that, by smoothness assumptions over $\mathbf{c}\mapsto \min \Arrowvert x-c_j\Arrowvert$,
we get, for any $\mathbf{c} \in (\mathbb{R}^d)^k$ and $\mathbf{c}^* \in \mathcal{M}$,               
               \[
               \gamma_\lambda(\mathbf{c},z)-\gamma_\lambda(\mathbf{c}^*,z)= \left\langle \mathbf{c}-\mathbf{c}^*),\nabla_{\c}\gamma_\lambda(\mathbf{c}^*,z) \right\rangle + \|\mathbf{c}-\mathbf{c}^*\|R_\lambda(\mathbf{c}^*,\mathbf{c}-\mathbf{c}^*,z),
               \]
               where, with \citet{Pollard82} we have:
$$
 \nabla_{\c}\gamma_\lambda(\mathbf{c}^*,z) = -2\left( \int\frac{1}{\lambda}\mathcal{K}_\eta\left(\frac{z-x}{\lambda}\right)(x-c^*_1)\mathbf{1}_{V^*_1}(x)dx,...,\int\frac{1}{\lambda}\mathcal{K}_\eta\left(\frac{z-x}{\lambda}\right)(x-c^*_k)\mathbf{1}_{V^*_k}(x)dx \right) $$
and $R_\lambda(          \mathbf{c}^*,\mathbf{c}-\mathbf{c}^*,z)$ such that:
$$
|R_\lambda(\mathbf{c}^*,\mathbf{c}-\mathbf{c}^*,z)|\leq\| \c - \c^* \|^{-1} \left(\left\langle\mathbf{c}-\mathbf{c}^*,\nabla_{\c}\gamma_\lambda(\mathbf{c}^*,z)\right\rangle+\max_{j=1,\ldots k}(|\|z-\mathbf{c}\|-\|x-\mathbf{c}^*\|\right).
$$
Splitting the expectation in two parts, we obtain:
\beqn
\label{dec}
&&\hspace{-0.9cm}\mathbb{E}  \sup_{\c^* \in \mathcal{M}, \|\mathbf{c}-\mathbf{c}^*\|^2 \leq \delta} {|\tilde{P}_n-\tilde{P}|(\gamma_\lambda(\mathbf{c}^*,.) - \gamma_\lambda(\mathbf{c},.))} \leq 
               \mathbb{E}   \sup_{\c^* \in \mathcal{M}, \|\mathbf{c}-\mathbf{c}^*\|^2 \leq \delta} { |\tilde{P}_n-\tilde{P}| \left\langle \mathbf{c}^*-\mathbf{c}, \nabla_{\c}\gamma_\lambda(\mathbf{c}^*,.) \right\rangle  }   \nonumber\\
               &+& \sqrt{\delta}  \mathbb{E} \sup_{\c^* \in \mathcal{M}, \|\mathbf{c}-\mathbf{c}^*\|^2 \leq \delta} {|\tilde{P}_n-\tilde{P}|(-R_\lambda(\mathbf{c}^*,\mathbf{c}-\mathbf{c}^*,. ))}    
\eeqn
 To bound the first term is this decomposition, consider the random variable
 \beqnn
 Z_n=(\tilde{P}_n-\tilde{P}) \left\langle \mathbf{c}^*-\mathbf{c}, \nabla_{\c}\gamma_\lambda(\mathbf{c}^*,.) \right\rangle =\frac{2}{n}\sum_{u=1}^k\sum_{j=1}^d(c_{u,j}-c^*_{u,j})\sum_{i=1}^n\int_{V_u}\frac{1}{\lambda}\mathcal{K}_\eta\left(\frac{Z_i-x}{\lambda}\right)(x_j-c_{u,j})dx.
 \eeqnn
By a simple Hoeffding's inequality, $Z_n$ is a subgaussian random variable. Its variance can be bounded as follows:
\beqnn
\mathrm{var} Z_n&=&\frac{4}{n}\sum_{u=1}^k\sum_{j=1}^d(c_{u,j}-c^*_{u,j})^2\mathrm{var}\int_{V_u}\frac{1}{\lambda}\mathcal{K}_\eta\left(\frac{Z-x}{\lambda}\right)(x_j-c_{u,j})dx\\
&\leq &\frac{4}{n}\delta \E\left(\int_{V_{u^+}}\frac{1}{\lambda}\mathcal{K}_\eta\left(\frac{Z-x}{\lambda}\right)(x_j-c_{{u^+},j})dx\right)^2\\
&\leq &C\frac{4}{n}\delta \int\left|\mathcal{F}\left[\frac{1}{\lambda}\mathcal{K}_\eta\left(\frac{\cdot}{\lambda}\right)\right](t)\right|^2\left|\mathcal{F}[(\pi_j-c_{{u^+},j})\mathrm{1}_{V_{u^+}}](t)\right|^2dt\\
&\leq &C\frac{4}{n}\delta\Pi_{i=1}^d\lambda_i^{-2\beta_i} \int_{V_{u^+}}(x_j-c_{u^+,j})^2dx\\
&\leq & C\Pi_{i=1}^d\lambda_i^{-2\beta_i}\frac{4}{n}\delta,
\eeqnn
where ${u^+}=\arg\max_u\int_{V_u}\frac{1}{\lambda}\mathcal{K}_\eta\left(\frac{Z-x}{\lambda}\right)(x_j-c_{u,j})dx$ and $\pi_j:x\mapsto x_j$, and where we use arguments originally stated in \citet{pinkfloyds} for compactly supported $\mathcal{F}[\mathcal{K}]$. We hence have using for instance a maximal inequality due to Massart \citet[Part 6.1]{saintflour}:
\beqnn
               \mathbb{E} \left (  \sup_{\c^* \in \mathcal{M}, \|\mathbf{c}-\mathbf{c}^*\|^2 \leq \delta} { (\tilde{P}_n-\tilde{P}) \left\langle \mathbf{c}^*-\mathbf{c}, \nabla_{\c}\gamma_\lambda(\mathbf{c}^*,.) \right\rangle  } \right ) 
                \leq
               C\frac{\Pi_{i=1}^d\lambda_i^{-\beta_i}}{\sqrt{n}}\sqrt{\delta}.
\eeqnn
                We obtain for the first term in (\ref{dec} the right order. To prove that the second term in (\ref{dec}) is smaller, note that from \citet{Pollard82}, we have:
\beqnn
|R_\lambda(\mathbf{c}^*,\mathbf{c}-\mathbf{c}^*,z)|&\leq&\| \c - \c^* \|^{-1} \left(\left\langle\mathbf{c}-\mathbf{c}^*,\nabla_{\c}\gamma_\lambda(\mathbf{c}^*,z)\right\rangle+\max_{j=1,\ldots k}(|\|z-\mathbf{c}_j\|^2-\|z-\mathbf{c}_j^*\|^2|\right)\\
&\leq &  \|\nabla_{\c}\gamma_\lambda(\mathbf{c}^*,z)\|+\| \c - \c^* \|^{-1}\sum_{j=1,\ldots k}|\|z-\mathbf{c}_j\|^2-\|z-\mathbf{c}_j^*\|^2|\\
&\leq &  C(\Pi_{i=1}^d\lambda_i^{-\beta_i}+\no z\no)
\eeqnn
we we use in last line:
$$
\|\nabla_{\c}\gamma_\lambda(\mathbf{c}^*,z)\|^2=4\sum_{j,k}\left(\int\frac{1}{\lambda}\mathcal{K}_\eta\left(\frac{z-x}{\lambda}\right)(x_j-c^*_{u,j})\mathbf{1}_{V^*_u}(x)dx\right)^2\leq C\Pi_{i=1}^d\lambda_i^{-2\beta_i}.
$$ 
Hence it is possible to apply a chaining argument as in \citet{levrard} to the class
$$
\mathcal{F}=\{R_\lambda(\mathbf{c}^*,\mathbf{c}-\mathbf{c}^*,\cdot),\mathbf{c}^*\in\mathcal{M},c\in\R^{kd}:\| \c-\c^*\|\leq\sqrt{\delta}\},
$$
which have an enveloppe function $F(\cdot)\leq C(\Pi_{i=1}^d\lambda_i^{-\beta_i}+\| \cdot\|)\in L_2(\tilde{P})$ provided that $\E\|\epsilon\|^2<\infty$.
\subsection{Proof of Theorem \ref{thm:application}}
Thr proof of Theorem \ref{thm:application} is divided into two steps. Using Theorem \ref{mainresult}, we can bound the excess risk when $|\mathcal{M}|=1$. For the general case of a finite numbers of optimal clusters $|\mathcal{M}|\geq 2$, we need to introduce a more sophisticated localization explain in \citet[Section 4]{kolt}.

\paragraph{First case: $|\mathcal{M}|=1$.\\}
The proof follows the proof of Theorem \ref{mainresult}. Using the previous notations, we have:
\beqnn
U_n^\lambda(\delta,t)\leq C\left(\phi_n^\lambda(\mathcal{C},\delta)+\sqrt{\frac{t}{n}\phi_n^\lambda(\mathcal{C},\delta)(1+K(\lambda))}+\sqrt{\frac{t}{n}}D^\lambda(\mathcal{C},\delta)+\frac{t}{n}\right).
\eeqnn
Using Lemma \ref{lemap}, gathering with Lemma \ref{complexity}, we have for $d=1$ for simplicity:
\beqnn
\phi_n^\lambda(\mathcal{C},\delta)&\leq &\E \sup_{\c,\c'\in\mathcal{C} (\delta)}|\tilde{P}_n-\tilde{P}|[\gamma_\lambda(\c)-\gamma_\lambda(\c')]\\
&\leq&\E \sup_{\|\mathbf{c}-\c^*\|^2\leq \delta}|\tilde{P}_n-\tilde{P}|[\gamma_\lambda(\c)-\gamma_\lambda(\c^*)]+\E \sup_{\|\mathbf{c}'-\c^*\|^2\leq \delta}|\tilde{P}_n-\tilde{P}|[\gamma_\lambda(\c')-\gamma_\lambda(\c^*)]\\
&\leq &C\frac{\lambda^{-\beta}}{\sqrt{n}}\delta^{\frac{1}{2}},
\eeqnn
where $\c^*$ is the unique minimizer of the clustering risk. Moreover, by uniqueness of $\c^*$, we can write from Lemma \ref{lemap}:
\beqnn
D^\lambda(\GG,\delta)&:=&\sup_{\c,\c'\in\mathcal{C}(\delta)}\sqrt{\tilde{P}(\gamma_\lambda(\c)-\gamma_\lambda(\c'))^2}\\
&\leq & C_1\lambda^{-\beta}\sup_{\c,\c'\in\mathcal{C}(\delta)}\sqrt{P(\gamma(\c)-\gamma(\c'))^2}\\
&\leq &C_1\lambda^{-\beta}\sup_{\c\in\mathcal{C}(\delta)}\sqrt{P(\gamma(\c)-\gamma(\c^*))^2}+\sup_{\c'\in\mathcal{C}(\delta)}\sqrt{P(\gamma(\c')-\gamma(\c^*))^2}\\
&\leq &2C_1\lambda^{-\beta}\sqrt{\delta}.
\eeqnn
It follows:
\beqnn
U_n^\lambda(\delta,t)\leq C\left(\frac{\lambda^{-\beta}}{\sqrt{n}}\delta^{\frac{1}{2}}+\frac{\lambda^{-\beta/2}}{n^{3/4}}\delta^{\frac{1}{4}}\sqrt{\lambda^{-\beta-1/2}t}+\sqrt{\frac{t}{n}} \lambda^{-\beta}\delta^{\frac{1}{2}}+\frac{t}{n}\right).
\eeqnn
We hence have the result applying Lemma \ref{interm} with the choice of $\lambda$ precised in Theorem \ref{thm:application}.
\paragraph{Second case: $|\mathcal{M}|\geq 2$\\}
When the infimum is not unique, the diameter $D^\lambda(\mathcal{G},\delta)$ does not necessary tend to zero when $\delta\to 0$. We hence introduce the more sophisticated geometric characteristic $r(\sigma,\delta)$ from \citet{kolt} defined as:
\beqnn
r(\sigma,\delta)=\sup_{\c\in\mathcal{C}(\delta)}\inf_{\c'\in\mathcal{C}(\sigma)}\sqrt{\tilde{P}(\gamma_\lambda(\c)-\gamma_\lambda(\c'))^2},\,\,0<\sigma\leq \delta.
\eeqnn
It is clear that $r(\sigma,\delta)\leq D^\lambda(\mathcal{G},\delta)$ and for $\delta\to 0$, we have $r(\sigma,\delta)\to 0$. The idea of the proof of Theorem \ref{thm:application} is to use a modified version of Theorem \ref{mainresult} using $r(\sigma,\delta)$ instead of $D^\lambda(\mathcal{G},\delta)$. Following \citet[Theorem 4]{kolt}, we can use a modified version of Lemma \ref{interm} in order to guarantee the upper bounds of Theorem \ref{mainresult} when $|\mathcal{M}|\geq 2$. To this end, we have to check for $d=1$ for simplicity:
\beqn
\label{newmodulus}
\lim_{\epsilon\to 0}\E\sup_{g\in\mathcal{G}(\sigma)}\sup_{g'\in\mathcal{G}(\delta):P(l(g)-l(g'))^2\leq r(\sigma,\delta)+\epsilon}\left|(\tilde{P}-\tilde{P}_n)(l_\lambda(g)-l_\lambda(g'))\right|\leq C \frac{\lambda^{-\beta}}{\sqrt{n}}\delta^{\frac{1-\rho}{2}},
\eeqn
and 
\beqn
\label{newdiam}
r(\sigma,\delta)\sqrt{\frac{t}{n}}\leq C\lambda^{-\beta}\sqrt{\frac{t\delta}{n}}.
\eeqn
Note that from Lemma \ref{lemap} and Lemma \ref{complexity}, it is clear that (\ref{newmodulus}) holds since:
\beqnn
&&\E\sup_{\c\in\mathcal{C}(\sigma)}\sup_{\c'\in\mathcal{C}(\delta):P(\gamma(\c)-\gamma(\c'))^2\leq r(\sigma,\delta)+\epsilon}\left|(\tilde{P}-\tilde{P}_n)(\gamma_\lambda(g)-\gamma_\lambda(g'))\right|\\
&&\leq \E\sup_{\c\in\mathcal{C}(\sigma),c^*\in\mathcal{M}}\left|(\tilde{P}-\tilde{P}_n)(\gamma_\lambda(\c)-\gamma_\lambda(\c^*))\right|+\E\sup_{\c'\in\mathcal{C}(\delta)}\left|(\tilde{P}-\tilde{P}_n)(\gamma_\lambda(\c')-\gamma_\lambda(\c^*(\c')))\right|\\
&&\leq 2\mathbb{E} \sup_{(\c,\c^*)\in\mathcal{C}\times\mathcal{M}, \|\mathbf{c}-\mathbf{c}^*\|^2 \leq c\delta} \left|(\tilde{P}_n-\tilde{P})(\gamma_\lambda(\mathbf{c}^*) - \gamma_\lambda(\mathbf{c}))\right| \\
&&\leq  C\lambda^{-\beta}\frac{\sqrt{\delta}}{\sqrt{n}}.
\eeqnn
Finally (\ref{newdiam}) holds since we have with Lemma \ref{lemap}, $\forall \c\in\mathcal{C}(\delta),\c'\in\mathcal{C}(\sigma)$:
\beqnn
\sqrt{\tilde{P}(\gamma_\lambda(\c)-\gamma_\lambda(\c'))^2}&\leq& C\lambda^{-\beta}\sqrt{P(\gamma(\c)-\gamma(\c')))^2}\\
&\leq&C\lambda^{-\beta}\left(\sqrt{P(\gamma(\c)-\gamma(\c^*(\c))^2}+\sqrt{P(\gamma(\c')-\gamma(\c^*(\c)))^2}\right)\\
&\leq & C\lambda^{-\beta}\sqrt{\delta},
\eeqnn
provided that $\sigma\leq \delta\leq 1.$
         \bibliographystyle{plain}
         \bibliography{referencejmlr}

\begin{thebibliography}{30}
\providecommand{\natexlab}[1]{#1}
\providecommand{\url}[1]{\texttt{#1}}
\expandafter\ifx\csname urlstyle\endcsname\relax
  \providecommand{\doi}[1]{doi: #1}\else
  \providecommand{\doi}{doi: \begingroup \urlstyle{rm}\Url}\fi

\bibitem[Antos et~al.(2005)Antos, Gy{\"o}rfi, and Gy{\"o}rgy]{gg}
A.~Antos, L.~Gy{\"o}rfi, and A.~Gy{\"o}rgy.
\newblock Individual convergence rates in empirical vector quantizer design.
\newblock \emph{IEEE Trans. Inform. Theory}, 51 (11), 2005.

\bibitem[Bartlett and Mendelson(2006)]{empimini}
P.L. Bartlett and S.~Mendelson.
\newblock Empirical minimization.
\newblock \emph{Probability Theory and Related Fields}, 135 (3):\penalty0
  311--334, 2006.

\bibitem[Bartlett et~al.(1998)Bartlett, Linder, and Lugosi]{lbclustering}
P.L. Bartlett, T.~Linder, and G.~Lugosi.
\newblock The minimax distortion redundancy in empirical quantizer design.
\newblock \emph{IEEE Trans. Inform. Theory}, 44 (5), 1998.

\bibitem[Bartlett et~al.(2005)Bartlett, Bousquet, and
  Mendelson]{localrademacher}
P.L. Bartlett, O.~Bousquet, and S.~Mendelson.
\newblock Local rademacher complexities.
\newblock \emph{The Annals of Statistics}, 33 (4):\penalty0 1497--1537, 2005.

\bibitem[Biau et~al.(2008)Biau, Devroye, and Lugosi]{biau}
G.~Biau, L.~Devroye, and G.~Lugosi.
\newblock On the performances of clustering in hilbert spaces.
\newblock \emph{IEEE Trans. Inform. Theory}, 54 (2), 2008.

\bibitem[Blanchard et~al.(2008)Blanchard, Bousquet, and Massart]{svm}
G.~Blanchard, O.~Bousquet, and P.~Massart.
\newblock Statistical performance of support vector machines.
\newblock \emph{The Annals of Statistics}, 36 (2):\penalty0 489--531, 2008.

\bibitem[Bousquet(2002)]{bousquet}
O.~Bousquet.
\newblock A bennet concentration inequality and its application to suprema of
  empirical processes.
\newblock \emph{C.R. Acad. SCI. Paris Ser. I Math}, 334:\penalty0 495--500,
  2002.

\bibitem[Butucea(2007)]{butucea}
C.~Butucea.
\newblock goodness-of-fit testing and quadratic functionnal estimation from
  indirect observations.
\newblock \emph{The Annals of Statistics}, 35:\penalty0 1907--1930, 2007.

\bibitem[Cavalier and Hengartner(2005)]{cavalierhengartner}
L.~Cavalier and N.W. Hengartner.
\newblock Adaptative estimation for inverse problems with noisy operators.
\newblock \emph{Inverse Problems}, 21 (4):\penalty0 1345--1361, 2005.

\bibitem[Fan(1991)]{Fan}
J.~Fan.
\newblock On the optimal rates of convergence for nonparametric deconvolution
  problems.
\newblock \emph{Annals of Statistics}, 19:\penalty0 1257--1272, 1991.

\bibitem[Graf and Luschgy(2000)]{existence}
Siegfried Graf and Harald Luschgy.
\newblock \emph{Foundation of quantization for probability distributions}.
\newblock Springer-Verlag, 2000.
\newblock Lecture Notes in Mathematics, volume 1730.

\bibitem[Hartigan(1975)]{hartigan75}
J.A. Hartigan.
\newblock \emph{Clustering algorithms}.
\newblock Wiley, 1975.

\bibitem[Koltchinskii(2006)]{kolt}
V.~Koltchinskii.
\newblock Local rademacher complexities and oracle inequalties in risk
  minimization.
\newblock \emph{The Annals of Statistics}, 34 (6):\penalty0 2593--2656, 2006.

\bibitem[Levrard(2012)]{levrard}
C.~Levrard.
\newblock Fast rates for empirical vector quantization.
\newblock \emph{hal.inria.fr/hal-00664068}, 2012.

\bibitem[Linder et~al.(1994)Linder, Lugosi, and Zeger]{llz}
T.~Linder, G.~Lugosi, and K.~Zeger.
\newblock Rates of convergence in the source coding theorem, in empirical
  quantizer design, and in universal lossy source coding.
\newblock \emph{IEEE Trans. Inform. Theory}, 40 (6), 1994.

\bibitem[Loustau(2011)]{noisykolt}
S.~Loustau.
\newblock Statistical learning with indirect observations.
\newblock \emph{http://hal.archives-ouvertes.fr/hal-00664125}, 2011.

\bibitem[Loustau and Marteau(2011)]{pinkfloyds}
S.~Loustau and C.~Marteau.
\newblock Discriminant analysis with errors in variables.
\newblock \emph{http://hal.archives-ouvertes.fr/hal-00660383}, 2011.

\bibitem[Mammen and Tsybakov(1999)]{mammen}
E.~Mammen and A.B. Tsybakov.
\newblock Smooth discrimination analysis.
\newblock \emph{The Annals of Statistics}, 27 (6):\penalty0 1808--1829, 1999.

\bibitem[Marteau(2006)]{marteau1}
C.~Marteau.
\newblock Regularization of inverse problems with unkown operator.
\newblock \emph{Mathematical Methods of Statistics}, 15 (4):\penalty0 415--443,
  2006.

\bibitem[Massart(2000)]{toulouse}
P.~Massart.
\newblock Some applications of concentration inequalities to statistics.
\newblock \emph{Ann. Fac. Sci. Toulouse Math.}, 9 (2):\penalty0 245--303, 2000.

\bibitem[Massart(2007)]{saintflour}
P.~Massart.
\newblock Concentration inequalities and model selection.
\newblock Ecole d'\'et\'e de Probabilit\'es de Saint-Flour 2003. Lecture Notes
  in Mathematics, Springer, 2007.

\bibitem[Massart and N\'ed\'elec(2006)]{nedelec}
P.~Massart and E.~N\'ed\'elec.
\newblock Risk bounds for statistical learning.
\newblock \emph{The Annals of Statistics}, 34 (5):\penalty0 2326--2366, 2006.

\bibitem[Meister(2009)]{meister}
A.~Meister.
\newblock \emph{Deconvolution problems in nonparametric statistics}.
\newblock Springer-Verlag, 2009.

\bibitem[Mendelson(2003)]{mendelsonkernel}
S.~Mendelson.
\newblock On the performance of kernel classes.
\newblock \emph{Journal of Machine Learning Research}, 4:\penalty0 759--771,
  2003.

\bibitem[Pollard(1981)]{pollard81}
D.~Pollard.
\newblock Strong consistency of k-means clustering.
\newblock \emph{The Annals of Statistics}, 9 (1), 1981.

\bibitem[Pollard(1982)]{Pollard82}
D.~Pollard.
\newblock A central limit theorem for $k$-means clustering.
\newblock \emph{The Annals of Probability}, 10 (4), 1982.

\bibitem[Tsybakov(2004)]{tsybakov2004}
A.B. Tsybakov.
\newblock Optimal aggregation of classifiers in statistical learning.
\newblock \emph{The Annals of Statistics}, 32 (1):\penalty0 135--166, 2004.

\bibitem[{Van De Geer}(2000)]{vdg}
S.~{Van De Geer}.
\newblock \emph{Empirical Processes in M-estimation}.
\newblock Cambridge University Press, 2000.

\bibitem[Vapnik(1982)]{vapnik82}
V.~Vapnik.
\newblock \emph{Estimation of Dependances Based on Empirical Data}.
\newblock Springer Verlag, 1982.

\bibitem[Vapnik(2000)]{vapnik2000}
V.~Vapnik.
\newblock \emph{The Nature of Statistical Learning Theory}.
\newblock Statistics for Engineering and Information Science, Springer, 2000.

\end{thebibliography}
\end{document}